# LIMIT THEOREMS FOR A CLASS OF IDENTICALLY DISTRIBUTED RANDOM VARIABLES

By Patrizia Berti, Luca Pratelli and Pietro Rigo

*Universita' di Modena e Reggio-Emilia, Accademia Navale di Livorno and Universita' di Pavia*

A new type of stochastic dependence for a sequence of random variables is introduced and studied. Precisely, $(X_n)_{n \geq 1}$ is said to be conditionally identically distributed (c.i.d.), with respect to a filtration $(\mathcal{G}_n)_{n \geq 0}$, if it is adapted to $(\mathcal{G}_n)_{n \geq 0}$ and, for each $n \geq 0$, $(X_k)_{k>n}$ is identically distributed given the past $\mathcal{G}_n$. In case $\mathcal{G}_0 = \{\varnothing, \Omega\}$ and $\mathcal{G}_n = \sigma(X_1, \ldots, X_n)$, a result of Kallenberg implies that $(X_n)_{n \geq 1}$ is exchangeable if and only if it is stationary and c.i.d. After giving some natural examples of nonexchangeable c.i.d. sequences, it is shown that $(X_n)_{n \geq 1}$ is exchangeable if and only if $(X_{\tau(n)})_{n \geq 1}$ is c.i.d. for any finite permutation $\tau$ of $\{1, 2, \ldots\}$, and that the distribution of a c.i.d. sequence agrees with an exchangeable law on a certain sub-$\sigma$-field. Moreover, $(1/n) \sum_{k=1}^n X_k$ converges a.s. and in $L^1$ whenever $(X_n)_{n \geq 1}$ is (real-valued) c.i.d. and $E[|X_1|] < \infty$. As to the CLT, three types of random centering are considered. One such centering, significant in Bayesian prediction and discrete time filtering, is $E[X_{n+1}|\mathcal{G}_n]$. For each centering, convergence in distribution of the corresponding empirical process is analyzed under uniform distance.

**1. Introduction and motivations.** In this paper a new type of stochastic dependence for a sequence $(X_n)_{n \geq 1}$ of random variables is introduced and studied. Precisely, suppose the $X_n$ are defined on the probability space $(\Omega, \mathcal{A}, P)$, take values in the measurable space $(E, \mathcal{E})$, and are adapted to a filtration $\mathcal{G} = (\mathcal{G}_n)_{n \geq 0}$. Then, $(X_n)_n$ is said to be *conditionally identically distributed with respect to* $\mathcal{G}$, abbreviated as $\mathcal{G}$-c.i.d., whenever

$$(1) \qquad E[f(X_k)|\mathcal{G}_n] = E[f(X_{n+1})|\mathcal{G}_n] \qquad \text{a.s.}$$

for all $k > n \geq 0$ and all bounded measurable $f : E \to \mathbb{R}$.

Roughly speaking, (1) means that, at each time $n \geq 0$, the future observations $(X_k)_{k>n}$ are identically distributed given the past $\mathcal{G}_n$. In case $\mathcal{G} = \mathcal{G}^X$, where $\mathcal{G}_0^X = \{\varnothing, \Omega\}$ and $\mathcal{G}_n^X = \sigma(X_1, \ldots, X_n)$, the filtration is not mentioned at all and $(X_n)_n$ is just called c.i.d. Clearly, if $(X_n)_n$ is $\mathcal{G}$-c.i.d., then it is c.i.d. and identically distributed.

Two obvious equivalent formulations of (1) are

$$(2) \qquad X_k \sim X_{n+1} \qquad \text{under } P(\cdot|H)$$

for all $k > n \geq 0$ and all events $H \in \mathcal{G}_n$ with $P(H) > 0$,

where "$\sim$" means "distributed as," and

$$(3) \qquad (E[f(X_{n+1})|\mathcal{G}_n])_{n \geq 0} \text{ is a } \mathcal{G}\text{-martingale}$$

for every bounded measurable $f : E \to \mathbb{R}$.

By general results on martingales, condition (3) can be written as $E[f(X_{T+1})] = E[f(X_1)]$ for all bounded measurable $f$ and all finite $\mathcal{G}$-stopping times $T$ (where $\mathcal{G}$-stopping times take values in $\{0, 1, \ldots, \infty\}$). Say that a $\mathcal{G}$-stopping time $S$ is *predictable* in case $S = T + 1$ for some $\mathcal{G}$-stopping time $T$. Then, one more equivalent formulation of (1) is

$$(4) \qquad X_S \sim X_1 \qquad \text{for each finite predictable } \mathcal{G}\text{-stopping time } S.$$

Note also that, when $\mathcal{G} = \mathcal{G}^X$, conditions (1)–(4) all reduce to

$$(5) \qquad [X_1, \ldots, X_n, X_{n+2}] \sim [X_1, \ldots, X_n, X_{n+1}] \qquad \text{for all } n \geq 0.$$

Exchangeable sequences meet (5) and, thus, are c.i.d. Indeed, exchangeability is the most significant case of conditional identity in distribution. C.i.d. sequences, however, need not be exchangeable. In fact, by a remarkable result of Kallenberg [(1988), Proposition 2.1], exchangeability amounts to stationarity and condition (5). In Kallenberg's paper (cf. Proposition 2.2), it is also shown that conditions (3)–(5) are equivalent in case $\mathcal{G} = \mathcal{G}^X$. However, apart from these results, condition (5) is not systematically investigated.

In the present paper, instead, we focus on $\mathcal{G}$-c.i.d. sequences. As a first motivation, we give some examples where conditional identity in distribution naturally arises while exchangeability may fail.

EXAMPLE 1.1 (Stopping and sampling). Let $X_n = Z_{T \wedge n}$, where $(Z_n)_n$ is exchangeable and $T$ is a random variable with values in $\{1, 2, \ldots, \infty\}$. Then, $(X_n)_n$ is not exchangeable apart from trivial cases, but it is c.i.d. under natural conditions on $T$. In fact, if $(Z_n)_n$ is c.i.d. (and not necessarily exchangeable), then $(X_n)_n$ is c.i.d. whenever $T$ is independent of $(Z_n)_n$, or whenever $T$ is a predictable stopping time for $\mathcal{G}^Z$. Thus, typically, conditional identity in distribution is preserved under stopping while exchangeability is not.

We now prove that $(X_n)_n$ is c.i.d. if $(Z_n)_n$ is c.i.d. and $T$ is a predictable stopping time for $\mathcal{G}^Z$. If $S$ is a finite predictable stopping time for $\mathcal{G}^Z$, then $T \wedge S$ is a finite predictable stopping time for $\mathcal{G}^Z$, and since $(Z_n)_n$ is c.i.d., one obtains

$$X_S = Z_{T \wedge S} \sim Z_1 = X_1.$$

Since $(X_n)_n$ is adapted to $\mathcal{G}^Z$, condition (4) implies that $(X_n)_n$ is $\mathcal{G}^Z$-c.i.d. and, in particular, it is c.i.d. Next, conditional identity in distribution is also preserved under (strictly increasing) sampling. That is, if $(Z_n)_n$ is c.i.d. and $T_1 < T_2 < \cdots$ are finite predictable stopping times for $\mathcal{G}^Z$, then $(X_n)_n = (Z_{T_n})$ is c.i.d. To prove the latter fact, fix a finite predictable stopping time $S$ for $\mathcal{G}^X$. Since $T_1 < T_2 < \cdots$, one has $\{T_j = n, S = j\} \in \mathcal{G}_{n-1}^Z$ for all $j, n \geq 1$, and this implies $\{T_S = n\} = \bigcup_{j=1}^{\infty} \{T_j = n, S = j\} \in \mathcal{G}_{n-1}^Z$. It follows that $T_S$ is a finite predictable stopping time for $\mathcal{G}^Z$, and since $(Z_n)_n$ is c.i.d., one obtains

$$X_S = Z_{T_S} \sim Z_1 \sim Z_{T_1} = X_1.$$

If $(X_n)_n$ is stationary and converges in probability, then $X_n = X_1$ a.s. for all $n$. In Example 1.1, if $T$ is a.s. finite, then $(X_n)_n = (Z_{T \wedge n})_n$ is definitively constant with probability 1 and, thus, it converges a.s. but in a trivial way. The next example exhibits a c.i.d. (nonexchangeable) sequence which converges a.s. in a nontrivial way.

EXAMPLE 1.2 (Compensated sum of independent random variables). Given the real numbers $0 < b_1 \leq b_2 \leq b_3 \leq \cdots < c$, let us define $\gamma_{ii} = c$ and $\gamma_{ij} = b_i \wedge b_j$ for $i \neq j$. On noting that $\Gamma_n = (\gamma_{ij})_{1 \leq i, j \leq n}$ is a symmetric positive definite matrix, $(X_n)_n$ can be taken such that $[X_1, \ldots, X_n] \sim \mathcal{N}(\mathbf{0}, \Gamma_n)$ for each $n \geq 1$. Then,

$$[X_1, \ldots, X_n, X_{n+2}] \sim \mathcal{N}(\mathbf{0}, \Gamma_{n+1}) \qquad \text{for all } n \geq 0,$$

that is, $(X_n)_n$ is c.i.d. [by condition (5)]. However,

$$E[(X_n - X_m)^2] = 2(c - b_n \wedge b_m).$$

Thus, $(X_n)_n$ is not stationary unless $b_n = b_1$ for all $n$, and $E[(X_n - X)^2] \to 0$ if $c = \lim_n b_n$, for some random variable $X$. Further, $X_n \to X$ a.s. whenever $\sum_n (c - b_n)^r < \infty$ for some $r > 0$ (since $E[|X_n - X|^{2r}] = \gamma_r (c - b_n)^r$ for some constant $\gamma_r$). To explain the title of the example, we note that it is a particular case of the following general scheme. Let $(Z_n)_n$, $(U_n)_n$ be independent sequences of independent real random variables and let

$$X_n = \sum_{i=1}^{n} Z_i + U_n, \qquad \mathcal{G}_n = \sigma(Z_1, U_1, \ldots, Z_n, U_n), \qquad \mathcal{G}_0 = \{\varnothing, \Omega\}.$$

Suppose also that $U_n$ compensates $\sum_{i=1}^{n} Z_i$, in the sense that $X_n \sim X_1$ for all $n$, and that the characteristic function $\phi_{X_1}$ of $X_1$ is null on a set with void interior. Fix $k > n \geq 0$ and a bounded Borel function $f : \mathbb{R} \to \mathbb{R}$. Then

$$E[f(X_k)|\mathcal{G}_n] = \int f\left(x + \sum_{i=1}^{n} Z_i\right)\mu_k(dx) \qquad \text{a.s.},$$

where $\mu_k$ is the distribution of $X_k - \sum_{i=1}^{n} Z_i$. But $\mu_k = \mu_{n+1}$, due to $\phi_{X_1}$ is null on a set with void interior and, thus, $(X_n)_n$ is $\mathcal{G}$-c.i.d. For instance, given any nondegenerate and infinitely divisible law $\mu$, the sequences $(Z_n)_n$ and $(U_n)_n$ can be taken such that the resulting $(X_n)_n$ is c.i.d., nonexchangeable with $X_1 \sim \mu$. Finally, to recover the first part of the example, just take $Z_n \sim \mathcal{N}(0, b_n - b_{n-1})$ and $U_n \sim \mathcal{N}(0, c - b_n)$, where $b_0 = 0$.

EXAMPLE 1.3 (Modified Pólya urns). An urn contains $w > 0$ white and $r > 0$ red balls. At each time $n \geq 1$, a ball is drawn and then replaced together with $d_n$ more balls of the same color. Let $X_n$ be the indicator of the event {white ball at time $n$}. Then $E(X_1) = w/(w + r)$ and

$$E[X_{n+1}|X_1, d_1, \ldots, X_n, d_n] = \frac{w + \sum_{i=1}^{n} d_i X_i}{w + r + \sum_{i=1}^{n} d_i} \qquad \text{a.s. for all } n \geq 1.$$

In the usual Pólya scheme, $d_n = d_1$ for all $n$, where $d_1 \geq 1$ is a fixed integer, and $(X_n)_n$ turns out to be exchangeable. Here, instead, we let $(d_n)_n$ be any sequence of random variables, with values in $\{1, 2, \ldots\}$, satisfying the following:

(i) $d_n$ is independent of $\sigma(X_i, d_j : i \leq n, j < n)$ for all $n \geq 1$, or
(ii) $d_1$ is degenerate and $\sigma(d_n) \subset \sigma(X_1, \ldots, X_{n-1})$ for all $n \geq 2$.

Then $(X_n)_n$ is c.i.d. but (apart from particular cases) nonexchangeable. For instance, if all the $d_n$ are degenerate, $(X_n)_n$ is not exchangeable unless $d_n = d_1$ for all $n$. To prove that $(X_n)_n$ is c.i.d., it is enough to check that $(E[X_{n+1}|\mathcal{G}_n])_{n \geq 0}$ is a $\mathcal{G}$-martingale for some filtration $\mathcal{G} \supset \mathcal{G}^X$. Suppose (i) holds and let

$$\mathcal{G}_0 = \{\varnothing, \Omega\}, \qquad \mathcal{G}_n = \sigma(X_1, d_1, \ldots, X_n, d_n, d_{n+1}) \qquad \text{for } n \geq 1.$$

For $n = 0$, (i) implies $E[E[X_2|\mathcal{G}_1]|\mathcal{G}_0] = E[X_2] = E[X_1] = E[X_1|\mathcal{G}_0]$ a.s. For $n \geq 1$, (i) gives $E[X_{n+1}|\mathcal{G}_n] = E[X_{n+1}|X_1, d_1, \ldots, X_n, d_n]$ a.s. Since $d_{n+1}$ is $\mathcal{G}_n$-measurable, it follows that

$$E[E[X_{n+2}|\mathcal{G}_{n+1}]|\mathcal{G}_n] = \frac{w + \sum_{i=1}^{n} d_i X_i + d_{n+1} E[X_{n+1}|\mathcal{G}_n]}{w + r + \sum_{i=1}^{n+1} d_i}$$

$$= E[X_{n+1}|\mathcal{G}_n] \qquad \text{a.s.}$$

A similar argument works under (ii), after setting $\mathcal{G} = \mathcal{G}^X$.

There is a second reason for studying $\mathcal{G}$-c.i.d. sequences, in addition to their possible utility in modelling real phenomena. Indeed, conditional identity in distribution is a basic assumption in *uniform limit theorems for predictive inference and empirical processes from dependent data*.

Precisely, suppose $E$ is a Polish space, $\mathcal{E} = \mathcal{B}(E)$ and $(X_n)_n$ is any sequence of random variables. Given a class $\mathcal{D}$ of bounded measurable functions on $E$, let

$$a_n(f) = E[f(X_{n+1})|\mathcal{G}_n] \qquad \text{for all } f \in \mathcal{D},$$

be the so-called *predictive measure*. In various problems, mainly in Bayesian predictive inference, discrete time filtering and sequential procedures, the main goal is just evaluating $a_n$, and good approximations $\tilde{a}_n$ for $a_n$ are needed. See, for instance, Algoet (1992, 1995), Ould-Said (1997), Modha and Masry (1998), Berti and Rigo (2002) and Berti, Mattei and Rigo (2002). Usually, $\tilde{a}_n$ is asked to meet a consistency condition of the type $\sup_{f \in \mathcal{D}} |\tilde{a}_n(f) - a_n(f)| \to 0$ a.s. A further request is that, for suitable normalizing constants $c_n$, the limiting distribution of $c_n(\tilde{a}_n - a_n)$ can be evaluated. Here, $c_n(\tilde{a}_n - a_n)$ is viewed as a process (indexed by $\mathcal{D}$) with paths in $l^\infty(\mathcal{D})$, the space of bounded functions on $\mathcal{D}$ equipped with uniform distance; see van der Vaart and Wellner (1996). In this framework, possible choices for $\tilde{a}_n$ and $c_n$ are the empirical measure $\mu_n = \frac{1}{n}\sum_{i=1}^n \delta_{X_i}$ and $c_n = \sqrt{n}$. So, it is of some interest to give conditions for

$$(6) \qquad \sup_{f \in \mathcal{D}} |\mu_n(f) - a_n(f)| \to 0 \qquad \text{a.s.}$$

$$(7) \qquad \sqrt{n}(\mu_n - a_n) \text{ converges in distribution to some known limit.}$$

Now, assuming that $(X_n)_n$ is $\mathcal{G}$-c.i.d. is fundamental for both (6) and (7). As to (6), we refer to Berti, Mattei and Rigo (2002). As to (7), one of the concerns of this paper is proving it for $\mathcal{G}$-c.i.d. sequences; see Section 4. Note also that (7) implies (6) if a.s. convergence is weakened to convergence in probability.

To sum up, conditional identity in distribution seems interesting enough to deserve a systematic study, both from the theoretical and the applied points of view. This task is accomplished here from the first point of view, with special attention to limit theorems.

The paper is organized in three sections. In Section 2, a few basic facts are listed. Among other things, a c.i.d. sequence meets a SLLN, is asymptotically exchangeable, and its probability distribution agrees with an exchangeable law on a certain sub-$\sigma$-field of $\mathcal{E}^\infty$. Moreover, $(X_n)_n$ is exchangeable if and only if $(X_{\tau(n)})_n$ is c.i.d. for any (finite) permutation $\tau$ of $\{1, 2, \dots\}$. Section 3 includes versions of the CLT for $\mathcal{G}$-c.i.d. sequences. Let $f: E \to \mathbb{R}$ be a measurable function. Stable convergence (in particular, convergence in distribution) of $\sqrt{n}(\frac{1}{n}\sum_{i=1}^n f(X_i) - L_n)$ is investigated for three different choices of the random centering $L_n$. In particular, conditions are given for convergence in distribution of

$$\sqrt{n}(\mu_n(f) - a_n(f)) = \sqrt{n}\left(\frac{1}{n}\sum_{i=1}^n f(X_i) - E[f(X_{n+1})|\mathcal{G}_n]\right).$$

Such conditions, incidentally, work in Examples 1.2 and 1.3. Section 4 is devoted to uniform limit theorems. For each centering considered in Section 3, convergence in distribution of the corresponding empirical process is investigated under uniform distance. General statements for $\mathcal{G}$-c.i.d. sequences are obtained which, among other things, yield interesting (and possibly new) results in the particular case of exchangeable sequences.

**2. Preliminary results and the SLLN.** Let $\mathcal{H}$ be the class of measurable functions $f: E \to \mathbb{R}$ such that $E[|f(X_1)|] < \infty$. Our starting point is the following simple lemma.

LEMMA 2.1. *Let $(X_n)_n$ be $\mathcal{G}$-c.i.d. Then, for each $f \in \mathcal{H}$, there is an integrable random variable $V_f$ such that $E[f(X_{n+1})|\mathcal{G}_n] \to V_f$, a.s. and in $L^1$, and*

$$(8) \qquad E[V_f|\mathcal{G}_n] = E[f(X_{n+1})|\mathcal{G}_n] \qquad \textit{a.s. for every } n \geq 0.$$

*Moreover, if $f_1, \ldots, f_k$ are bounded elements of $\mathcal{H}$, $k > 1$, then*

$$(9) \qquad E\left[\prod_{j=1}^k f_j(X_{n+j}) \Big| \mathcal{G}_n\right] \to \prod_{j=1}^k V_{f_j} \qquad \textit{a.s. and in } L^1.$$

PROOF. By (3), $(E[f(X_{n+1})|\mathcal{G}_n])_{n\geq 0}$ is a $\mathcal{G}$-martingale, and it is uniformly integrable since the $X_n$ are identically distributed. Hence, $E[f(X_{n+1})|\mathcal{G}_n] \to V_f$, a.s. and in $L^1$, for some random variable $V_f$. In particular, $V_f$ closes the martingale $(E[f(X_{n+1})|\mathcal{G}_n])_{n\geq 0}$ and, thus, condition (8) holds. As to (9), since $f_1, \ldots, f_k$ are bounded, it is enough to show a.s. convergence. Arguing by induction, suppose that

$$E\left[\prod_{j=1}^{k-1} f_j(X_{n+j}) \Big| \mathcal{G}_n\right] \to \prod_{j=1}^{k-1} V_{f_j} \qquad \text{a.s.}$$

Let $D_n = E[(V_{f_k} - E[V_{f_k}|\mathcal{G}_n])\prod_{j=1}^{k-1} f_j(X_{n+j})|\mathcal{G}_n]$. Since $f_1, \ldots, f_k$ are bounded and $E[V_{f_k}|\mathcal{G}_n] \to V_{f_k}$ a.s., it follows that $D_n \to 0$ a.s. Hence, (8) and the inductive assumption imply

$$E\left[\prod_{j=1}^k f_j(X_{n+j}) \Big| \mathcal{G}_n\right]$$

$$= E\left[\prod_{j=1}^{k-1} f_j(X_{n+j}) E[f_k(X_{n+k})|\mathcal{G}_{n+k-1}]\Big|\mathcal{G}_n\right]$$

$$= E\left[\prod_{j=1}^{k-1} f_j(X_{n+j}) V_{f_k}\Big|\mathcal{G}_n\right]$$

$$= D_n + E[V_{f_k}|\mathcal{G}_n] E\left[\prod_{j=1}^{k-1} f_j(X_{n+j})\Big|\mathcal{G}_n\right] \to \prod_{j=1}^{k} V_{f_j} \qquad \text{a.s.} \qquad \square$$

Among other things, Lemma 2.1 has implications as regards convergence in $\sigma(L^1, L^\infty)$ of c.i.d. sequences. Recall that, for real integrable random variables $Y_n$ and $Y$ on the same probability space, $Y_n \to Y$ in $\sigma(L^1, L^\infty)$ means $E[Y_n Z] \to E[Y Z]$ for each bounded random variable $Z$. Then, $f(X_n) \to V_f$ in $\sigma(L^1, L^\infty)$ whenever $(X_n)_n$ is $\mathcal{G}$-c.i.d. and $f \in \mathcal{H}$. Fix, in fact, a bounded random variable $Z$. By standard arguments, for proving $E[Zf(X_n)] \to E[ZV_f]$ it can be assumed that $Z$ is $\mathcal{G}_m$-measurable for some $m$, and in this case Lemma 2.1 yields

$$E[ZV_f] = \lim_n E[ZE[f(X_n)|\mathcal{G}_{n-1}]] = \lim_n E[Zf(X_n)].$$

Moreover, in exactly the same way, Lemma 2.1 also implies that

$$\prod_{j=1}^{k} f_j(X_{n+j}) \to \prod_{j=1}^{k} V_{f_j} \qquad \text{in } \sigma(L^1, L^\infty)$$

whenever $(X_n)_n$ is $\mathcal{G}$-c.i.d. and $f_1, \ldots, f_k$ are bounded elements of $\mathcal{H}$. From now on, when $(X_n)_n$ is c.i.d. and $f \in \mathcal{H}$, $V_f$ always denotes a version of the limit in $\sigma(L^1, L^\infty)$ of $(f(X_n))_n$.

If $(X_n)_n$ is $\mathcal{G}$-c.i.d., then $(f(X_n))_n$ is still $\mathcal{G}$-c.i.d. for each measurable $f$ on $E$, while $(g(X_n, X_{n+1}, \ldots))_n$ can fail to be c.i.d. if $g$ is measurable on $E^\infty$; see, for instance, Example 1.2. Nevertheless, $(g(X_n, X_{n+1}, \ldots))_n$ obeys a SLLN for various choices of $g$, for instance, for $g$ of the form $g(x) = \prod_{j=1}^{k} f_j(x_j)$, $x = (x_1, x_2, \ldots) \in E^\infty$.

THEOREM 2.2 (SLLN). Let $(X_n)_n$ be c.i.d. If $f_1, \ldots, f_k \in \mathcal{H}$ and $f_1, \ldots, f_{k-1}$ are bounded, $k \geq 1$, then

$$(10) \qquad \frac{1}{n} \sum_{i=0}^{n-1} \prod_{j=1}^{k} f_j(X_{i+j}) \to \prod_{j=1}^{k} V_{f_j} \qquad \text{a.s. and in } L^1.$$

In particular, $\frac{1}{n} \sum_{i=1}^{n} f(X_i) \to V_f$, a.s. and in $L^1$, whenever $f \in \mathcal{H}$.

PROOF. Let $U_i = \prod_{j=1}^{k} f_j(X_{i+j})$, $i \geq 0$. Since $(U_i)_i$ is uniformly integrable, it is enough to prove a.s. convergence, and, to this end, it can be assumed $f_j \geq 0$ for all $j$. To begin with, suppose also that $f_k$ is bounded, and let

$$Z_n = \sum_{i=0}^{n-1} \frac{U_i - E[U_i|\mathcal{G}_{i+k-1}^X]}{i+1}.$$

Then, $(Z_n)_n$ is a martingale with respect to $(\mathcal{G}_{n+k-1}^X)_n$, and since $f_1, \ldots, f_k$ are all bounded, one has $\sup_n E[Z_n^2] < \infty$. Hence, $(Z_n)_n$ converges a.s., and an application of the Kronecker lemma gives

(11) $$\frac{1}{n} \sum_{i=0}^{n-1} (U_i - E[U_i|\mathcal{G}_{i+k-1}^X]) \to 0 \qquad \text{a.s.}$$

For $k = 1$, one has $\frac{1}{n} \sum_{i=0}^{n-1} E[U_i|\mathcal{G}_{i+k-1}^X] = \frac{1}{n} \sum_{i=0}^{n-1} E[f_1(X_{i+1})|\mathcal{G}_i^X] \to V_{f_1}$ a.s. by Lemma 2.1 and, thus, (11) implies (10). Arguing by induction, suppose that

(12) $$\frac{1}{n} \sum_{i=0}^{n-1} \prod_{j=1}^{k-1} f_j(X_{i+j}) \to \prod_{j=1}^{k-1} V_{f_j} \qquad \text{a.s.}$$

Recall that, if $(a_n)_n$ and $(b_n)_n$ are any real sequences, then $\frac{1}{n} \sum_{i=0}^{n-1} a_i b_i \to ab$ whenever $\frac{1}{n} \sum_{i=0}^{n-1} a_i \to a$, $b_n \to b$ and $a_i \geq 0$ for all $i$. Hence,

$$\frac{1}{n} \sum_{i=0}^{n-1} E[U_i|\mathcal{G}_{i+k-1}^X]$$

$$= \frac{1}{n} \sum_{i=0}^{n-1} \prod_{j=1}^{k-1} f_j(X_{i+j}) E[f_k(X_{i+k})|\mathcal{G}_{i+k-1}^X] \to \prod_{j=1}^{k} V_{f_j} \qquad \text{a.s.}$$

Once again, (10) follows from (11), and this concludes the proof in the particular case where $f_k$ is bounded. If $f_k$ is not bounded, define $f_{k,m} = f_k I_{\{f_k \leq m\}}$. Then, $V_{f_{k,m}} \uparrow V_{f_k}$ a.s. as $m \to \infty$. Further, for each fixed $m_0$, one obtains

$$V_{f_k} \geq \limsup_m E[f_{k,m}(X_{m+k})|\mathcal{G}_{m+k-1}^X]$$

$$\geq \liminf_m E[f_{k,m}(X_{m+k})|\mathcal{G}_{m+k-1}^X]$$

$$\geq \liminf_m E[f_{k,m_0}(X_{m+k})|\mathcal{G}_{m+k-1}^X] = V_{f_{k,m_0}} \qquad \text{a.s.}$$

Thus, $V_{f_k} = \lim_m E[f_{k,m}(X_{m+k})|\mathcal{G}_{m+k-1}^X]$ a.s. Let $Y_i = f_{k,i}(X_{i+k}) \times \prod_{j=1}^{k-1} f_j(X_{i+j})$. Since $f_1, \ldots, f_{k-1}$ are bounded, condition (12) holds by the first part of this proof. Therefore,

$$\frac{1}{n}\sum_{i=0}^{n-1}E[Y_i|\mathcal{G}_{i+k-1}^X] = \frac{1}{n}\sum_{i=0}^{n-1}\prod_{j=1}^{k-1}f_j(X_{i+j})E[f_{k,i}(X_{i+k})|\mathcal{G}_{i+k-1}^X]$$

$$\rightarrow \prod_{j=1}^{k}V_{f_j} \qquad \text{a.s.}$$

Further, $\sum_{i=0}^{\infty}P(U_i\neq Y_i)\leq\sum_{i=0}^{\infty}P(f_k(X_1)>i)\leq 1+E[f_k(X_1)]$ which implies $P(U_i\neq Y_i \text{ i.o.})=0$. Hence, it suffices showing that $\frac{1}{n}\sum_{i=0}^{n-1}(Y_i-E[Y_i|\mathcal{G}_{i+k-1}^X])\rightarrow 0$ a.s. In its turn, this follows from the Kronecker lemma, after noting that $H_n=\sum_{i=0}^{n-1}\frac{Y_i-E[Y_i|\mathcal{G}_{i+k-1}^X]}{i+1}$ is a martingale such that $\sup_n E[H_n^2]<\infty$. In fact, letting $a=\prod_{j=1}^{k-1}\sup|f_j|$, one obtains

$$\sup_n E[H_n^2]\leq\sum_{n=1}^{\infty}n^{-2}E[Y_{n-1}^2]$$

$$\leq 2a^2\sum_{n=2}^{\infty}n^{-2}\int_0^{n-1}xP(f_k(X_1)>x)\,dx$$

$$\leq 2a^2\sum_{n=2}^{\infty}n^{-2}\sum_{i=1}^{n-1}iP(f_k(X_1)>i-1)$$

$$=2a^2\sum_{i=1}^{\infty}iP(f_k(X_1)>i-1)\sum_{n=i+1}^{\infty}n^{-2}$$

$$\leq 2a^2\sum_{i=1}^{\infty}P(f_k(X_1)>i-1)$$

$$\leq 2a^2(1+E[f_k(X_1)]). \qquad \square$$

REMARK 2.3. Suppose $E$ is a Polish space, $\mathcal{E}=\mathcal{B}(E)$ and $(X_n)_n$ is c.i.d. By using Theorem 2.2 and Lemma 2.4, it is not hard to see that $\frac{1}{n}\sum_{i=1}^{n}g(X_i,X_{i+1},\dots)$ converges a.s. for each bounded continuous function $g$ on $E^{\infty}$. This result generally fails if $g$ is a bounded, Borel but not continuous function on $E^{\infty}$.

The remaining part of this section investigates to what extent conditional identity in distribution is connected with exchangeability. To this end, we collect here some notation and terminology from Aldous (1985). Given a Polish space $S$, let $C_b(S)$ denote the space of bounded continuous functions on $S$, $\mathbb{P}$ the set of probability measures on $\mathcal{B}(S)$, and $\Sigma$ the $\sigma$-field on $\mathbb{P}$ generated by the evaluation maps $p\mapsto p(B)$, for $B$ varying in $\mathcal{B}(S)$. A *random* *measure on $S$ is a measurable function $\gamma : (\Omega, \mathcal{A}) \to (\mathbb{P}, \Sigma)$.* Let $(Z_n)_n$ be a sequence of $S$-valued random variables on $(\Omega, \mathcal{A}, P)$. Say that $(Z_n)_n$ converges *stably* if, for every $H \in \mathcal{A}$ with $P(H) > 0$, $(Z_n)_n$ converges in distribution under $P(\cdot|H)$ to some law $\mu_H$. In this case, there is a random measure $\gamma$ on $S$ which represents each limit law $\mu_H$ as $\mu_H(\cdot) = \int \gamma(\omega)(\cdot) P(d\omega|H)$, and $(Z_n)_n$ is said *to converge stably with representing measure $\gamma$.* See also [Letta and Pratelli (1996)](). We recall that, if $(Z_n)_n$ is exchangeable, there is a random measure $\gamma$ on $S$ such that the product random measure

$$\gamma^\infty = \gamma \times \gamma \times \cdots$$

is a version of the conditional distribution of $(Z_n)_n$, given $\sigma(\gamma)$. Such $\gamma$ is called the *directing measure of $(Z_n)_n$.* For our purposes, a last simple fact should be stressed.

LEMMA 2.4. *Let $S$ be a Polish space and $(Z_n)_n$ any sequence of $S$-valued random variables on $(\Omega, \mathcal{A}, P)$. If $\sigma(Z_1, Z_2, \dots) \subset \sigma(\bigcup_n \mathcal{G}_n)$ and $(E[f(Z_{n+1})| \mathcal{G}_n])_n$ converges a.s. for each $f \in C_b(S)$, there is a random measure $\gamma$ on $S$ such that $(Z_n)_n$ converges stably with representing measure $\gamma$ and*

$$E[f(Z_{n+1})|\mathcal{G}_n] \to \int f(x) \gamma(\cdot)(dx) \qquad \text{a.s. for each } f \in C_b(S).$$

PROOF. Fix $H \in \mathcal{G}_m$ with $P(H) > 0$, $m \geq 1$. Then, $E_{P(\cdot|H)}[f(Z_n)]$ converges to a finite limit for all $f \in C_b(S)$, so that $(Z_n)_n$ converges in distribution under $P(\cdot|H)$. Since $\sigma(Z_1, Z_2, \dots) \subset \sigma(\bigcup_n \mathcal{G}_n)$, it follows that $(Z_n)_n$ converges in distribution under $P(\cdot|H)$ for each $H \in \mathcal{A}$ with $P(H) > 0$. Hence, there is a random measure $\gamma$ on $S$ such that $(Z_n)_n$ converges stably with representing measure $\gamma$. Such $\gamma$ can be taken $\sigma(Z_1, Z_2, \dots)$-measurable. Let $f \in C_b(S)$ and $L_f = \lim_n E[f(Z_{n+1})|\mathcal{G}_n]$ a.s. Then, $E[I_H L_f] = E[I_H \int f(x)\gamma(\cdot)(dx)]$ for all $H \in \bigcup_n \mathcal{G}_n$. Since $L_f$ and $\int f(x)\gamma(\cdot)(dx)$ are measurable with respect to $\sigma(\bigcup_n \mathcal{G}_n)$, one obtains $L_f = \int f(x)\gamma(\cdot)(dx)$ a.s. □

Suppose $E$ is a Polish space, $\mathcal{E} = \mathcal{B}(E)$ and $(X_n)_n$ is $\mathcal{G}$-c.i.d. In view of Lemmas [2.1]() and [2.4](), there is a random measure $\alpha$ on $E$ such that $V_f = \int f(t)\alpha(\cdot)(dt)$ a.s. for all $f \in C_b(E)$. By a monotone class argument, it follows that

$$V_{I_B}(\omega) = \alpha(\omega)(B) \qquad \text{for almost all } \omega$$

whenever $B \in \mathcal{E}$. In the sequel, $\alpha$ is called the *directing measure of $(X_n)_n$.* Such terminology, which is typical of exchangeable sequences, is motivated by at least two facts. First, by Theorem [2.2](), $\alpha$ is the a.s. weak limit of empirical measures,

$$\frac{1}{n}\sum_{i=1}^n \delta_{X_i(\omega)} \to \alpha(\omega) \qquad \text{weakly, for almost all } \omega.$$

Second, $(X_n)_n$ is asymptotically exchangeable and the exchangeable limit law is, in a sense, directed by $\alpha$. In fact, something more is true in the following:

THEOREM 2.5. *Suppose $E$ is a Polish space, $\mathcal{E} = \mathcal{B}(E)$ and $(X_n)_n$ is $\mathcal{G}$-c.i.d. Then, $[X_n, X_{n+1}, \ldots]$ converges stably with representing measure $\alpha^\infty$, where $\alpha$ is the directing measure of $(X_n)_n$, and*

$$
\begin{aligned}
(13) \quad & E[g(X_{n+1}, X_{n+2}, \ldots)|\mathcal{G}_n] \\
& \rightarrow \int g(x)\alpha^\infty(\cdot)(dx) \qquad \text{a.s. for each } g \in C_b(E^\infty).
\end{aligned}
$$

PROOF. For each $n \geq 0$, fix a regular version $\nu_n$ of the conditional distribution of $[X_{n+1}, X_{n+2}, \ldots]$ given $\mathcal{G}_n$, and define

$$
\nu_{n,k}(\omega)(B) = \nu_n(\omega)\{x \in E^\infty : x_{k-n} \in B\} \qquad \text{for all } k > n, \omega \in \Omega \text{ and } B \in \mathcal{E}.
$$

Let $\mu_n = \nu_{n,n+1}$. Since $(X_n)_n$ is $\mathcal{G}$-c.i.d. and $\nu_{n,k}$ is a version of the conditional distribution of $X_k$ given $\mathcal{G}_n$, one obtains $\nu_{n,k} = \mu_n$ a.s. for all $k > n \geq 0$. By Lemma 2.1, $\int f(t)\mu_n(\cdot)(dt) \rightarrow \int f(t)\alpha(\cdot)(dt)$ a.s. for each $f \in C_b(E)$, and this implies $\mu_n(\omega) \rightarrow \alpha(\omega)$ weakly for almost all $\omega$. Let $H_1 \in \mathcal{A}$ be such that $P(H_1) = 1$ and, for all $\omega \in H_1$, $(\mu_n(\omega))_n$ is tight and $\nu_{n,k}(\omega) = \mu_n(\omega)$ for all $k > n \geq 0$. Then, $(\nu_n(\omega))_n$ is tight, too, for all $\omega \in H_1$. Let $D_0 \subset C_b(E)$ be a countable convergence determining class for $E$, and let $D$ be the class of those functions $g$ on $E^\infty$ of the form $g(x) = \prod_{j=1}^k f_j(x_j)$, $x = (x_1, x_2, \ldots) \in E^\infty$, where $k \geq 1$ and $f_1, \ldots, f_k \in D_0$. For $g \in D$, say $g(x) = \prod_{j=1}^k f_j(x_j)$, Lemma 2.1 gives

$$
\int g(x)\alpha^\infty(\omega)(dx) = \prod_{j=1}^k \int f_j(t)\alpha(\omega)(dt) = \prod_{j=1}^k V_{f_j}(\omega) = \lim_n \int g(x)\nu_n(\omega)(dx)
$$

for almost all $\omega$. Since $D$ is countable, there is $H_2 \in \mathcal{A}$ with $P(H_2) = 1$ and $\int g(x)\alpha^\infty(\omega)(dx) = \lim_n \int g(x)\nu_n(\omega)(dx)$ for all $g \in D$ and $\omega \in H_2$. It follows that $\nu_n(\omega) \rightarrow \alpha^\infty(\omega)$ weakly for all $\omega \in H_1 \cap H_2$. This proves (13). To conclude the proof, it suffices to apply Lemma 2.4 with $S = E^\infty$ and $Z_n = [X_n, X_{n+1}, \ldots]$. □

By Theorem 2.5, $[X_n, X_{n+1}, \ldots]$ converges in distribution under $P(\cdot|H)$ to the exchangeable law $\mu_H(\cdot) = \int \alpha^\infty(\omega)(\cdot)P(d\omega|H)$ whenever $(X_n)_n$ is c.i.d., $H \in \mathcal{A}$ and $P(H) > 0$. An alternative proof of this fact could be given by results of Aldous (1985). Incidentally, we also note that Theorem 2.5 directly implies Kallenberg's result that $(X_n)_n$ is exchangeable if and only if it is stationary and c.i.d.; see Section 1. In fact, if $(X_n)_n$ is stationary and c.i.d., the distribution of $[X_n, X_{n+1}, \ldots]$ does not depend on $n$ (by stationarity) and converges weakly to the exchangeable law $\mu_\Omega$ [by Theorem 2.5, since $(X_n)_n$ is c.i.d.].

One more consequence of Theorems 2.2 and 2.5 is that any c.i.d. law on $\mathcal{E}^\infty$ is exchangeable on a suitable sub-$\sigma$-field. Let $\pi_n$ be the $n$th coordinate projection on $E^\infty$ and $\mathcal{V}$ the $\sigma$-field on $E^\infty$ generated by

$$\limsup_n \frac{1}{n}[f(\pi_1) + \cdots + f(\pi_n)] \qquad \text{for all bounded } f \in \mathcal{H}.$$

THEOREM 2.6. *Suppose $E$ is a Polish space, $\mathcal{E} = \mathcal{B}(E)$ and $(X_n)_n$ is c.i.d. Then, the probability distribution $\lambda$ of $(X_n)_n$ coincides on $\mathcal{V}$ with the exchangeable law $\mu(\cdot) = \int \alpha^\infty(\omega)(\cdot) P(d\omega)$, where $\alpha$ is the directing measure of $(X_n)_n$.*

PROOF. By Theorem 2.5, $[X_n, X_{n+1}, \ldots]$ converges in distribution to $\mu$. Given $f_1, \ldots, f_k \in C_b(E)$, this fact and Lemma 2.1 imply

$$\int \prod_{j=1}^k f_j \circ \pi_j \, d\mu = \lim_n E\left[\prod_{j=1}^k f_j(X_{n+j})\right] = E\left[\prod_{j=1}^k V_{f_j}\right].$$

Let $\mathcal{V}_0$ be the $\sigma$-field on $E^\infty$ generated by $\limsup_n \frac{1}{n}[f(\pi_1) + \cdots + f(\pi_n)]$ for all $f \in C_b(E)$, and let $h$ be any product of generators of $\mathcal{V}_0$, that is,

$$h = \prod_{j=1}^k \limsup_n \frac{1}{n} \sum_{i=1}^n f_j \circ \pi_i \qquad \text{where } f_1, \ldots, f_k \in C_b(E).$$

By Theorem 2.2, $\limsup_n \frac{1}{n} \sum_{i=1}^n f_j(X_i) = V_{f_j}$ a.s. and, thus, $\int h \, d\lambda = E[\prod_{j=1}^k V_{f_j}]$. On the other hand, exchangeability of $\mu$ implies

$$\int h \, d\mu = \int \lim_n \prod_{j=1}^k \frac{1}{n} \sum_{i=1}^n f_j \circ \pi_i \, d\mu = \lim_n \frac{1}{n^k} \int \prod_{j=1}^k \sum_{i=1}^n f_j \circ \pi_i \, d\mu$$

$$= \int \prod_{j=1}^k f_j \circ \pi_j \, d\mu = E\left[\prod_{j=1}^k V_{f_j}\right] = \int h \, d\lambda.$$

Hence, $\lambda = \mu$ on $\mathcal{V}_0$. To conclude the proof, it is sufficient showing that $\mathcal{V} \subset \sigma(\mathcal{V}_0 \cup \mathcal{N})$, where $\mathcal{N} = \{A \in \mathcal{E}^\infty : \lambda(A) = \mu(A) = 0\}$. Let $\nu = (\lambda + \mu)/2$. Given a bounded measurable $\phi$ on $E$ and $\varepsilon > 0$, there is $f \in C_b(E)$ such that $\int |\phi(\pi_1) - f(\pi_1)| \, d\nu < \varepsilon$. Since $(\pi_n)_n$ is c.i.d. under $\nu$, Theorem 2.2 implies

$$\int \left| \limsup_n \frac{1}{n} \sum_{i=1}^n \phi(\pi_i) - \limsup_n \frac{1}{n} \sum_{i=1}^n f(\pi_i) \right| d\nu$$

$$= \lim_n \int \left| \frac{1}{n} \sum_{i=1}^n (\phi(\pi_i) - f(\pi_i)) \right| d\nu$$

$$\leq \int |\phi(\pi_1) - f(\pi_1)| \, d\nu < \varepsilon.$$

This concludes the proof. $\quad\square$

REMARK 2.7. In Theorem 2.6, $\mathcal{V}$ cannot be replaced by the shift-invariant $\sigma$-field of $(\pi_n)_n$. In fact, if the distribution of $(X_n)_n$ agrees with an exchangeable law $\mu$ on the shift-invariant $\sigma$-field of $(\pi_n)_n$, then

$$P(\exists \lim X_n) = \mu(\exists \lim \pi_n)$$

$$= \mu(\exists m \text{ with } \pi_n = \pi_m \text{ for all } n \geq m)$$

$$= P(\exists m \text{ with } X_n = X_m \text{ for all } n \geq m),$$

where the second equality is due to exchangeability of $\mu$. But, there are c.i.d. sequences for which $P(\exists \lim X_n) = 1 > 0 = P(\exists m \text{ with } X_n = X_m \text{ for all } n \geq m)$, for instance, the one exhibited in Example 1.2. It follows that, unlike in the exchangeable case, the shift-invariant $\sigma$-field of $(\pi_n)_n$ and $\mathcal{V}$ have not the same completion under an arbitrary c.i.d. law on $\mathcal{E}^\infty$.

We close this section with a characterization of exchangeability in terms of conditional identity in distribution.

THEOREM 2.8. *Let $\mathcal{I}$ be the shift-invariant $\sigma$-field of $(X_n)_n$. The following statements are equivalent:*

(i) *$(X_n)_n$ is exchangeable;*
(ii) *For any $H \in \mathcal{I}$ with $P(H) > 0$, $(X_n)_n$ is c.i.d. under $P(\cdot|H)$;*
(iii) *For any (finite) permutation $\tau$ of $\{1, 2, \dots\}$, $(X_{\tau(n)})_n$ is c.i.d.*

PROOF. (i) $\Rightarrow$ (ii). Obvious. (ii) $\Rightarrow$ (iii). Fix $f \in \mathcal{H}$ and note that, by Theorem 2.2, $V_f$ can be taken $\mathcal{I}$-measurable. Hence, by (ii), $E[f(X_1)|\mathcal{I}] = V_f$ a.s. Further, let $n \geq 1$, $H \in \mathcal{G}_n^X$ and $K \in \mathcal{I}$. For all $k > n$, condition (ii) implies

$$E[f(X_k)I_H I_K] = E[f(X_{n+1})I_H I_K].$$

Since $f(X_k) \to V_f$ in $\sigma(L^1, L^\infty)$ as $k \to \infty$, one also obtains $E[V_f I_H I_K] = E[f(X_{n+1})I_H I_K]$. Since $V_f$ is $\mathcal{I}$-measurable, this implies

$$E[f(X_{n+1})|\mathcal{I}] = V_f = E[f(X_{n+1})|\mathcal{I} \vee \mathcal{G}_n^X] \qquad \text{a.s.}$$

Thus, $(X_n)_n$ is exchangeable and, clearly, this implies condition (iii).

(iii) $\Rightarrow$ (i). We prove that, for any $n_1, \ldots, n_p$ distinct integers, any $r \geq 1$ and any $m > \max(n_1, \ldots, n_p)$,

$$(14) \qquad E\left[\prod_{j=1}^{r} f_j(X_{m+j}) \prod_{l=1}^{p} g_l(X_{n_l})\right] = E\left[\prod_{j=1}^{r} f_j(X_{m+j}) \prod_{l=1}^{p} g_l(X_l)\right],$$

where $f_1, \ldots, f_r, g_1, \ldots, g_p$ are bounded elements of $\mathcal{H}$. We argue by induction on $p \geq 1$. When $p = 1$, condition (14) follows from applying (5) to the c.i.d. sequence $(X_{\tau(n)})_n$, where $\tau$ is a finite permutation such that $\tau(j) = m + j$ for $j = 1, \ldots, r$, $\tau(r+1) = n_1$ and $\tau(r+2) = 1$. Suppose now that (14) holds for some $p$. We have to prove

$$(15) \qquad E\left[\prod_{j=1}^{k} f_j(X_{m+j}) \prod_{l=1}^{p+1} g_l(X_{n_l})\right] = E\left[\prod_{j=1}^{k} f_j(X_{m+j}) \prod_{l=1}^{p+1} g_l(X_l)\right]$$

for any $k \geq 1$ and any $m > \max(n_1, \ldots, n_p, n_{p+1})$. Let $\tau_1$ denote a finite permutation such that $\tau_1(j) = m + j$ for $j = 1, \ldots, k$, $\tau_1(k+l) = n_l$ for $l = 1, \ldots, p+1$ and $\tau_1(k+2+p) = m+k+1$. Since $(X_{\tau_1(n)})_n$ is c.i.d., one has

$$E\left[\prod_{j=1}^{k} f_j(X_{m+j}) \prod_{l=1}^{p+1} g_l(X_{n_l})\right] = E\left[\prod_{j=1}^{k} f_j(X_{m+j}) \prod_{l=1}^{p} g_l(X_{n_l}) g_{p+1}(X_{n_{p+1}})\right]$$

$$= E\left[\prod_{j=1}^{k} f_j(X_{m+j}) \prod_{l=1}^{p} g_l(X_{n_l}) g_{p+1}(X_{m+k+1})\right].$$

Let $\tau_2$ be a finite permutation such that $\tau_2(j) = m + j$ for $j = 1, \ldots, k$, $\tau_2(k+l) = l$ for $l = 1, \ldots, p+1$ and $\tau_2(k+2+p) = m+k+1$. Since $(X_{\tau_2(n)})_n$ is c.i.d., one also has

$$E\left[\prod_{j=1}^{k} f_j(X_{m+j}) \prod_{l=1}^{p+1} g_l(X_l)\right]$$

$$= E\left[\prod_{j=1}^{k} f_j(X_{m+j}) \prod_{l=1}^{p} g_l(X_l) g_{p+1}(X_{p+1})\right]$$

$$= E\left[\prod_{j=1}^{k} f_j(X_{m+j}) \prod_{l=1}^{p} g_l(X_l) g_{p+1}(X_{m+k+1})\right].$$

Hence, (15) follows from (14) with $r = k+1$ and $f_{k+1} = g_{p+1}$. $\square$

## 3. Some CLTs for c.i.d. sequences.

In this section stable convergence (in particular, convergence in distribution) of

$$\sqrt{n}\left(\frac{1}{n}\sum_{k=1}^{n} f(X_k) - L_n\right)$$

is investigated for three different choices of the random centering $L_n$.

In all cases, our main tool is the following version of the martingale CLT; see Hall and Heyde (1980), Theorem 3.2, page 58. Let $\{Y_{nk} : n \geq 1, k = 1, \ldots, k_n\}$ be an array of real square integrable random variables, where $k_n \uparrow \infty$, and for all $n$, let $\mathcal{F}_{n0} \subset \mathcal{F}_{n1} \subset \cdots \subset \mathcal{F}_{nk_n} \subset \mathcal{A}$ be $\sigma$-fields with $\mathcal{F}_{n0} = \{\varnothing, \Omega\}$. If:

   (i) $\sigma(Y_{nk}) \subset \mathcal{F}_{nk}, E[Y_{nk}|\mathcal{F}_{n,k-1}] = 0$ a.s., $\mathcal{F}_{nk} \subset \mathcal{F}_{n+1,k}$,
   (ii) $\max_{1 \leq k \leq k_n} |Y_{nk}| \to 0$ in probability, $\sup_n E[\max_{1 \leq k \leq k_n} Y_{nk}^2] < \infty$,
   (iii) $\sum_{k=1}^{k_n} Y_{nk}^2 \to L$ in probability, for some real random variable $L$,

then $\sum_{k=1}^{k_n} Y_{nk}$ converges stably. Precisely, let $\mathcal{N}(0, c)$ denote the Gaussian law with mean 0 and variance $c \geq 0$, where $\mathcal{N}(0, 0) = \delta_0$. Then, $\sum_{k=1}^{k_n} Y_{nk}$ converges stably with representing measure $\mathcal{N}(0, L)$, that is, for each $H \in \mathcal{A}$ with $P(H) > 0$, $\sum_{k=1}^{k_n} Y_{nk}$ converges in distribution under $P(\cdot|H)$ to the probability law

$$\mu_H(\cdot) = \int \mathcal{N}(0, L(\omega))(\cdot) P(d\omega|H).$$

Let us start with the case $L_n = V_f$.

THEOREM 3.1 (CLT, case I).   *Suppose $(X_n)_n$ is c.i.d., $f$ and $f^2$ are in $\mathcal{H}$ and there exists an integer $m \geq 0$ such that $(f(X_{n+m}) - V_f)_n$ is c.i.d. Then*

$$W_{n,f} = \frac{1}{\sqrt{n}}(f(X_1) + \cdots + f(X_n) - nV_f)$$

*converges stably with representing measure $\mathcal{N}(0, V_{f^2} - (V_f)^2)$.*

PROOF.   For each $n > m$, define

$$Y_{nk} = n^{-1/2}(f(X_{k+m}) - V_f)$$

and $\mathcal{F}_{nk} = \sigma(Y_{n1}, \ldots, Y_{nk})$ for $k = 1, \ldots, n - m$. Since $(f(X_{n+m}) - V_f)_n$ is c.i.d., condition (i) holds. By Theorem 2.2, condition (iii) holds with $L = V_{f^2} - (V_f)^2$. As to (ii), first note that it can be equivalently written as

$$n^{-1/2} \max_{m < k \leq n} |f(X_k)| \to 0 \qquad \text{in probability}$$

and

$$\sup_n n^{-1} E\left[ \max_{m < k \leq n} f^2(X_k) \right] < \infty.$$

Fix $\varepsilon > 0$ and put $A_{nk} = \{|f(X_k)| > \varepsilon\sqrt{n}\,\}$. Then

$$P\left(\max_{m < k \leq n}|f(X_k)| > \varepsilon\sqrt{n}\right) \leq \frac{1}{\varepsilon^2 n}\sum_{k=m+1}^{n}E[I_{A_{nk}}f^2(X_k)]$$

$$= \frac{n-m}{\varepsilon^2 n}E[I_{A_{n1}}f^2(X_1)] \to 0$$

and

$$\sup_n \frac{1}{n}E\left[\max_{m < k \leq n}f^2(X_k)\right] \leq \sup_n \frac{1}{n}\sum_{k=m+1}^{n}E[f^2(X_k)] \leq E[f^2(X_1)]. \qquad \square$$

If $(X_n)_n$ is exchangeable, then $(f(X_{n+m}) - V_f)_n$ is exchangeable for all $m$, so that Theorem 3.1 applies. Generally, however, the assumption that $(f(X_{n+m}) - V_f)_n$ is c.i.d. for some $m$ cannot be dropped.

EXAMPLE 3.2 (Example 1.2 continued).   Let $X_n = \sum_{k=1}^{n}Z_k + U_n$, where $(Z_n)_n$ and $(U_n)_n$ are independent sequences of independent random variables, $Z_n \sim \mathcal{N}(0, b_n - b_{n-1})$, $U_n \sim \mathcal{N}(0, c - b_n)$, with $b_0 = 0$ and $b_n \uparrow c$. Let $f$ be the identity mapping, $\mathcal{G}_0 = \{\varnothing, \Omega\}$ and $\mathcal{G}_n = \sigma(Z_1, U_1, \ldots, Z_n, U_n)$. Then $(X_n)_n$ is $\mathcal{G}$-c.i.d. Since $X_k - U_k = \sum_{i=1}^{k}Z_i \to \sum_{i=1}^{\infty}Z_i$ a.s., Theorem 2.2 implies

$$V_f = \lim_n \frac{1}{n}\sum_{k=1}^{n}X_k = \lim_n \frac{1}{n}\sum_{k=1}^{n}(X_k - U_k) = \sum_{k=1}^{\infty}Z_k \qquad \text{a.s.}$$

Then

$$\sqrt{n}W_{n,f} = \sum_{k=1}^{n}((n-k+1)Z_k + U_k) - n\sum_{k=1}^{\infty}Z_k = \sum_{k=1}^{n}(U_k - (k-1)Z_k) - n\sum_{k>n}Z_k,$$

so that $W_{n,f}$ is Gaussian with

$$\mathrm{Var}[W_{n,f}] = \frac{1}{n}\sum_{k=1}^{n}(c - b_k + (k-1)^2(b_k - b_{k-1})) + n(c - b_n).$$

Hence, if $\lim_n n(c - b_n) = \infty$, then $\lim_n \mathrm{Var}[W_{n,f}] = \infty$. In that case, since $W_{n,f}$ is Gaussian for all $n$, $(W_{n,f})_n$ does not converge in distribution.

Let us turn to the second type of random centering, that is, $L_n = \frac{1}{n} \times \sum_{k=1}^{n}E[f(X_k)|\mathcal{G}_{k-1}]$. This is perhaps the less interesting of our choices of $L_n$, at least from the point of view of applications. Nevertheless, there are situations where such a choice of $L_n$ plays a role, for instance, in stochastic approximation, calibration and gambling; see Hanson and Russo (1981,

1986), Dawid (1982) and Berti and Rigo (2002). In any case, the following result is available (CLT, case II; we omit the straightforward proof). Let

$$B_{n,f} = \frac{1}{\sqrt{n}} \sum_{k=1}^{n} (f(X_k) - E[f(X_k)|\mathcal{G}_{k-1}]) \qquad \text{for all } f \in \mathcal{H}.$$

If $(X_n)_n$ is $\mathcal{G}$-c.i.d. and $f$ and $f^2$ are in $\mathcal{H}$, then

$$(B_{n,f})_n \text{ converges stably with representing measure } \mathcal{N}(0, V_{f^2} - (V_f)^2). \tag{16}$$

Finally, we consider the case $L_n = E[f(X_{n+1})|\mathcal{G}_n]$. From the point of view of statistical applications, mainly in Bayesian forecasting and discrete time filtering, this is perhaps the most significant case; see Section 1. Denote

$$C_{n,f} = \frac{1}{\sqrt{n}} (f(X_1) + \cdots + f(X_n) - nE[f(X_{n+1})|\mathcal{G}_n]) \qquad \text{for all } f \in \mathcal{H}.$$

THEOREM 3.3 (CLT, case III). *Suppose $(X_n)_n$ is $\mathcal{G}$-c.i.d., $f$ and $f^2$ are in $\mathcal{H}$ and $\sup_n E[C_{n,f}^2] < \infty$. If*

$$M_n = \frac{1}{n} \sum_{k=1}^{n} (f(X_k) - kE[f(X_{k+1})|\mathcal{G}_k] + (k-1)E[f(X_k)|\mathcal{G}_{k-1}])^2 \to \sigma^2 \qquad a.s.$$

*for some real random variable $\sigma^2$, then $(C_{n,f})_n$ converges stably with representing measure $\mathcal{N}(0, \sigma^2)$. Moreover, if*

$$\frac{1}{n} \sum_{k=1}^{n} k^2 (E[f(X_{k+1})|\mathcal{G}_k] - E[f(X_k)|\mathcal{G}_{k-1}])^2 \to 0 \qquad \text{in probability,}$$

*then $B_{n,f} - C_{n,f} \to 0$ in probability, and $(C_{n,f})_n$ converges stably with representing measure $\mathcal{N}(0, V_{f^2} - (V_f)^2)$.*

PROOF. Suppose first that $M_n \to \sigma^2$ a.s. For $n \geq 1$ and $k = 1, \ldots, n$, define $Y_{nk} = E[W_{n,f}|\mathcal{G}_k] - E[W_{n,f}|\mathcal{G}_{k-1}]$ and $\mathcal{F}_{nk} = \mathcal{G}_k$. Then condition (i) trivially holds, and since $E[W_{n,f}|\mathcal{G}_0] = 0$ a.s., one has $C_{n,f} = E[W_{n,f}|\mathcal{G}_n] = \sum_{k=1}^{n} Y_{nk}$. Hence, it is enough to check (ii) and (iii) with $L = \sigma^2$. On noting that

$$\sqrt{n} Y_{n,k} = f(X_k) - kE[f(X_{k+1})|\mathcal{G}_k] + (k-1)E[f(X_k)|\mathcal{G}_{k-1}],$$

one obtains $\sum_{k=1}^{n} Y_{nk}^2 = M_n \to \sigma^2$ a.s. Since $Y_{nn}^2 = M_n - \frac{n-1}{n} M_{n-1} \to 0$ a.s., it follows that $\max_{k \leq n} |Y_{nk}| \to 0$ a.s. Moreover,

$$E\left[\max_{k \leq n} Y_{nk}^2\right] \leq \sum_{k=1}^{n} E[(E[W_{n,f}|\mathcal{G}_k] - E[W_{n,f}|\mathcal{G}_{k-1}])^2]$$

$$= \sum_{k=1}^{n} E[E[W_{n,f}|\mathcal{G}_k]^2 - E[W_{n,f}|\mathcal{G}_{k-1}]^2]$$

$$= E[E[W_{n,f}|\mathcal{G}_n]^2] = E[C_{n,f}^2],$$

so that $(\max_{k \leq n} |Y_{nk}|)_n$ is bounded in $L^2$. Hence, (ii) and (iii) hold with $L = \sigma^2$, and this concludes the proof of the first part of the theorem. Next, to prove the second part, define $D_k = E[f(X_{k+1})|\mathcal{G}_k] - E[f(X_k)|\mathcal{G}_{k-1}]$ and suppose that $\frac{1}{n}\sum_{k=1}^{n} k^2 D_k^2 \to 0$ in probability. By (16), it is sufficient to see that $B_{n,f} - C_{n,f} \to 0$ in probability, and a direct calculation shows that $B_{n,f} - C_{n,f} = \frac{1}{\sqrt{n}}\sum_{k=1}^{n} k D_k$. For $n \geq 1$ and $k = 1, \ldots, n$, define $Y_{nk} = \frac{1}{\sqrt{n}} k D_k$ and $\mathcal{F}_{nk} = \mathcal{G}_k$. Then, (i) and (iii) hold with $L = 0$. In particular, $\{\max_{k \leq n} |Y_{nk}|\}^2 \leq \sum_{k=1}^{n} Y_{nk}^2 \to 0$ in probability. Moreover,

$$\max_{k \leq n} |Y_{nk}| \leq \frac{1}{\sqrt{n}} \max_{k \leq n} |f(X_k) - E[f(X_k)|\mathcal{G}_{k-1}] - k D_k|$$

$$+ \frac{1}{\sqrt{n}} \max_{k \leq n} |f(X_k) - E[f(X_k)|\mathcal{G}_{k-1}]|,$$

and both terms in the right-hand side are bounded in $L^2$. (Boundedness of the first term has been shown in the first part of the proof.) Hence, condition (ii) holds, and this implies $\frac{1}{\sqrt{n}}\sum_{k=1}^{n} k D_k = \sum_{k=1}^{n} Y_{nk} \to 0$ in probability. □

The assumption that $(C_{n,f})_n$ is bounded in $L^2$ surely holds if $(W_{n,f})_n$ is bounded in $L^2$ and, in turn, this is true if $(f(X_{n+m}) - V_f)_n$ is c.i.d. for some $m$. In particular, Theorem 3.3 implies that $(C_{n,f})_n$ converges stably whenever $(X_n)_n$ is exchangeable, $E[f(X_1)^2] < \infty$ and $(M_n)_n$ converges a.s. Here, it is tempting to conjecture that $(M_n)_n$ always converges a.s. in the exchangeable case, but we do not know whether this is true.

We close this section by applying the previous results to some of the examples in Section 1. Example 3.4 shows that, for c.i.d. nonexchangeable sequences, $(W_{n,f})_n$, $(B_{n,f})_n$ and $(C_{n,f})_n$ can have quite different asymptotic behaviors. Example 3.5 deals with modified Pólya urns, in the particular case where the extra balls $d_n$ are i.i.d.

EXAMPLE 3.4 (Example 1.2 continued). Let $X_n = \sum_{k=1}^{n} Z_k + U_n$, where $(Z_n)_n$ and $(U_n)_n$ are independent sequences of independent random variables, $Z_n \sim \mathcal{N}(0, b_n - b_{n-1})$, $U_n \sim \mathcal{N}(0, c - b_n)$, with $b_0 = 0$ and $b_n \uparrow c$. Let $f$ be the identity mapping, $\mathcal{G}_0 = \{\varnothing, \Omega\}$ and $\mathcal{G}_n = \sigma(Z_1, U_1, \ldots, Z_n, U_n)$. Then, as noted in Example 3.2, $(X_n)_n$ is $\mathcal{G}$-c.i.d., $V_f = \sum_{k=1}^{\infty} Z_k$ a.s. and $W_{n,f}$ is Gaussian with mean 0 and $\text{Var}[W_{n,f}] = \frac{1}{n}\sum_{k=1}^{n}(c - b_k + (k-1)^2(b_k - b_{k-1})) + n(c - b_n)$. Suppose now that

$$\sup_n n^2(b_n - b_{n-1}) < \infty \quad \text{and} \quad n(c - b_n) \to u$$

for some $u$. Then, a direct calculation shows that $\frac{1}{n}\sum_{k=1}^{n}(k-1)^2(b_k - b_{k-1}) \to u$, and $M_n = \frac{1}{n}\sum_{k=1}^{n}(U_k-(k-1)Z_k)^2 \to u$ a.s. Thus, $(W_{n,f})_n$ converges in distribution to $\mathcal{N}(0,2u)$ and, by Theorem 3.3, $(C_{n,f})_n$ converges in distribution to $\mathcal{N}(0,u)$. Finally,

$$V_{f^2} = \lim_n \frac{1}{n}\sum_{k=1}^{n} X_k^2 = \lim_n \frac{1}{n}\sum_{k=1}^{n}\left(\sum_{i=1}^{k} Z_i\right)^2 = \left(\sum_{i=1}^{\infty} Z_i\right)^2 = V_f^2 \qquad \text{a.s.}$$

and thus (16) yields $B_{n,f} \to 0$ in probability.

EXAMPLE 3.5 (Example 1.3 continued). Let $(X_n)_n$ and $(d_n)_n$ be as in Example 1.3, and let $\mathcal{G}_0 = \{\varnothing, \Omega\}$ and $\mathcal{G}_n = \sigma(X_1, d_1, \ldots, X_n, d_n, d_{n+1})$ for $n \geq 1$. Suppose that $d_n$ is independent of $\sigma(X_i, d_j : i \leq n, j < n)$ for all $n \geq 1$ [i.e., condition (i) holds] and the $d_n$ are identically distributed with $E[d_1^2] < \infty$. As shown in Example 1.3, $(X_n)_n$ is $\mathcal{G}$-c.i.d. Let $f$ be the identity mapping. By standard but long calculations, it can be shown that $(C_{n,f})_n$ is bounded in $L^2$ and $M_n \to \delta(V - V^2)$ a.s., where $\delta = \mathrm{Var}[d_1]/E[d_1]^2$ and $V = \lim_n \frac{1}{n}\sum_{k=1}^{n} X_k$ a.s. We refer to Berti, Pratelli and Rigo (2002) for details on such calculations. In any case, by Theorem 3.3, $(C_{n,f})_n$ converges stably with representing measure $\mathcal{N}(0, \delta(V - V^2))$.

**4. Uniform limit theorems.** In Section 3, given a $\mathcal{G}$-c.i.d. sequence $(X_n)_n$, convergence in distribution of

$$W_{n,f} = \frac{1}{\sqrt{n}}(f(X_1) + \cdots + f(X_n) - nV_f),$$

$$B_{n,f} = \frac{1}{\sqrt{n}}\sum_{k=1}^{n}(f(X_k) - E[f(X_k)|\mathcal{G}_{k-1}]),$$

$$C_{n,f} = \frac{1}{\sqrt{n}}(f(X_1) + \cdots + f(X_n) - nE[f(X_{n+1})|\mathcal{G}_n])$$

has been investigated for a *fixed* function $f$. In this section $W_n := \{W_{n,f} : f \in \mathcal{D}\}$, $B_n := \{B_{n,f} : f \in \mathcal{D}\}$ and $C_n := \{C_{n,f} : f \in \mathcal{D}\}$ are seen as processes, indexed by some class $\mathcal{D} \subset \mathcal{H}$ of functions, and their convergence in distribution is analyzed in the path space under uniform distance. Note that $W_n$, $B_n$ and $C_n$ all reduce to the usual empirical process whenever $X_n$ is independent of $\mathcal{G}_{n-1}$ for all $n$.

For the sake of simplicity, we do not deal with a general Donsker-class $\mathcal{D}$, but we focus on the particular case where

$$(E, \mathcal{E}) = (\mathbb{R}, \mathcal{B}(\mathbb{R})) \quad \text{and} \quad \mathcal{D} = \{I_{(-\infty, t]} : t \in \mathbb{R}\}.$$

Hence, letting

$$\mathcal{X} = \left\{ x : x \text{ is a real cadlag function on } \mathbb{R} \text{ and } \lim_{|t| \to \infty} x(t) = 0 \right\},$$

$$W_{n,t} = W_{n,I_{(-\infty,t]}}, \qquad B_{n,t} = B_{n,I_{(-\infty,t]}}, \qquad C_{n,t} = C_{n,I_{(-\infty,t]}},$$

the paths of $W_n$, $B_n$ and $C_n$ belong to $\mathcal{X}$ (up to modifications on $P$-null sets).

Throughout, $\mathcal{X}$ is equipped with uniform distance. We refer to the theory of weak convergence developed by Hoffmann-Jørgensen, van der Vaart and Wellner; see van der Vaart and Wellner (1996). Let $(\Omega', \mathcal{A}', P')$ be a probability space and $Z : \Omega' \to \mathcal{X}$ a random element of $\mathcal{X}$. Say that $Z$ is *measurable* if $\{Z \in B\} \in \mathcal{A}'$ for all Borel sets $B \subset \mathcal{X}$, and that $Z$ is *tight* if $Z$ is indistinguishable from a measurable random element with a tight probability distribution. If $Z$ and $Z'$ are both measurable and tight, $Z \sim Z'$ if and only if they have the same finite-dimensional distributions. If $Z$ is measurable and the $Z_n$ are arbitrary random elements of $\mathcal{X}$, then $Z_n \to Z$ in distribution means $E^*[f(Z_n)] \to E[f(Z)]$ for all $f \in C_b(\mathcal{X})$, where $E^*$ denotes outer expectation. If $Z$ is not measurable, but indistinguishable from a measurable random element $Z'$, then $Z_n \to Z$ in distribution stands for $Z_n \to Z'$ in distribution. Suppose the $Z_n$ are random processes on $(\Omega, \mathcal{A}, P)$ such that $(Z_{n,t_1}, \ldots, Z_{n,t_r})$ converges in distribution for all $t_1, \ldots, t_r \in \mathbb{R}$. Then, for $(Z_n)_n$ to converge in distribution to a tight limit, it is sufficient that, for all $\varepsilon, \eta > 0$, there is a finite partition $I_1, \ldots, I_m$ of $\mathbb{R}$ by right-open intervals such that

$$(17) \qquad \limsup_n P\left( \max_k \sup_{s,t \in I_k} |Z_{n,s} - Z_{n,t}| > \varepsilon \right) < \eta;$$

see van der Vaart and Wellner (1996), Theorems 1.5.4 and 1.5.6.

When $(X_n)_n$ is $\mathcal{G}$-c.i.d., a possible limit in distribution for $W_n$, $B_n$ and $C_n$ is a tight process whose distribution $\nu$ is given by

$$\nu\{ x \in \mathcal{X} : (x(t_1), \ldots, x(t_r)) \in A \} = \int \mathcal{N}(0, \Sigma(t_1, \ldots, t_r))(A)\, dP$$

for all $t_1, \ldots, t_r \in \mathbb{R}$ and $A \in \mathcal{B}(\mathbb{R}^r)$, where $\Sigma(t_1, \ldots, t_r)$ is a random covariance matrix. One significant particular case is the following. Let $\mathbb{G}^F$ denote a process, on some probability space, of the form

$$\mathbb{G}_t^F = \mathbb{G}_{F(t)}^0, \qquad \text{for all } t \in \mathbb{R},$$

where $\mathbb{G}^0$ is a standard Brownian bridge on $[0,1]$ and $F$ a random distribution function, independent of $\mathbb{G}^0$, satisfying

$$(F(t_1), \ldots, F(t_r)) \sim (\alpha(-\infty, t_1], \ldots, \alpha(-\infty, t_r]) \qquad \text{for all } t_1, \ldots, t_r \in \mathbb{R},$$

$\alpha$ being the directing measure of $(X_n)_n$. Then, $\mathbb{G}^F$ has finite-dimensional distributions of the type of $\nu$ with

$$\Sigma(t_1, \ldots, t_r) = (F(t_i \wedge t_j)(1 - F(t_i \vee t_j)))_{1 \leq i, j \leq r}.$$

Generally, $\mathbb{G}^F$ can fail to be measurable. However, $\mathbb{G}^F$ is measurable and tight whenever all the $F$-paths are continuous on $A^c$ for some fixed countable set $A \subset \mathbb{R}$.

Before stating results, we will give a technical lemma that is needed later on. It is presumably well known, and we provide a proof just to make the paper self-contained. Let us denote $\|x\| = \sup_t |x(t)|$ for all $x \in \mathcal{X}$.

LEMMA 4.1. *Let $Z$ be a tight random process with paths in $\mathcal{X}$ and $E[\|Z\|] < \infty$. Then, for all $\varepsilon > 0$, there is a finite partition $I_1, \ldots, I_m$ of $\mathbb{R}$ by right-open intervals such that*

$$E\left(\max_k \sup_{s, t \in I_k} |Z_s - Z_t|\right) < \varepsilon.$$

PROOF. It can be assumed that $Z$ is measurable. By tightness of $Z$ and integrability of $\|Z\|$, there is a compact $K$ such that $E[I_{\{Z \notin K\}}\|Z\|] < \varepsilon/5$. Let $x_1, \ldots, x_N \in K$ be such that $K \subset \bigcup_{i=1}^N B_i$, where $B_i$ is the ball with center $x_i$ and radius $\varepsilon/5$. Take a partition $I_1, \ldots, I_m$ of $\mathbb{R}$ by right-open intervals such that $\max_k \sup_{s, t \in I_k} |x_i(s) - x_i(t)| < \varepsilon/5$ for all $i = 1, \ldots, N$. Then

$$I_{\{Z \in K\}} \max_k \sup_{s, t \in I_k} |Z_s - Z_t| < (3/5)\varepsilon$$

and thus $E(\max_k \sup_{s, t \in I_k} |Z_s - Z_t|) \leq (3/5)\varepsilon + 2E[I_{\{Z \notin K\}}\|Z\|] < \varepsilon$. $\square$

Next, based on the results in Section 3, we give conditions for convergence in distribution of $(B_n)_n$ and $(C_n)_n$.

THEOREM 4.2. *If $(X_n)_n$ is $\mathcal{G}$-c.i.d. and $(B_n)_n$ meets condition (17), then $B_n \to \mathbb{G}^F$ in distribution and $\mathbb{G}^F$ is tight.*

PROOF. First note that $\mathbb{G}^F$ is measurable when $\mathcal{X}$ is equipped with the ball $\sigma$-field $\mathcal{U}$. Suppose the finite-dimensional distributions of $(B_n)_n$ converge weakly to those of $\mathbb{G}^F$. Then, since $(B_n)_n$ meets (17), $B_n \to Z$ in distribution for some measurable tight process $Z$ with the same finite-dimensional distributions of $\mathbb{G}^F$. Since $Z$ is tight, $Z \in A$ a.s. for some separable Borel set $A \subset \mathcal{X}$. Since $A \in \mathcal{U}$ (by separability) and the distributions of $\mathbb{G}^F$ and $Z$ agree on $\mathcal{U}$, one obtains $\mathbb{G}^F \in A$ a.s. Let $L = I_{\{\mathbb{G}^F \in A\}} F + I_{\{\mathbb{G}^F \notin A\}} H$, where $H$ is any fixed distribution function. Then $L$ is a random distribution function indistinguishable from $F$, and $\mathbb{G}^L$ is measurable and tight due to $\mathbb{G}^L$ having separable range and $\mathcal{X}$ is complete. Since $\mathbb{G}^F$ is indistinguishable from $\mathbb{G}^L$, it follows that $\mathbb{G}^F$ is tight, $Z \sim \mathbb{G}^L$ and $B_n \to \mathbb{G}^F$ in distribution. It remains to prove that the finite-dimensional distributions of $(B_n)_n$ converge weakly to those of $\mathbb{G}^F$. Fix $t_1, \dots, t_r, a_1, \dots, a_r \in \mathbb{R}$, define $f = \sum_{i=1}^{r} a_i I_{(-\infty, t_i]}$ and note that $\sum_{i=1}^{r} a_i \mathbb{G}^F_{t_i}$ has distribution

$$\mu(A) = \int \mathcal{N}(0, V_{f^2} - (V_f)^2)(A) \, dP, \qquad A \in \mathcal{B}(\mathbb{R}).$$

Hence, (16) implies

$$\sum_{i=1}^{r} a_i B_{n,t_i} = B_{n,f} \to \sum_{i=1}^{r} a_i \mathbb{G}^F_{t_i} \qquad \text{in distribution.}$$

By letting $a_1, \dots, a_r$ vary, one obtains $(B_{n,t_1}, \dots, B_{n,t_r}) \to (\mathbb{G}^F_{t_1}, \dots, \mathbb{G}^F_{t_r})$ in distribution. $\quad\square$

Convergence in distribution of $(C_n)_n$ needs more conditions. Furthermore, as suggested by Theorem 3.3, Examples 3.4 and 3.5, it may be that $C_n \to C$ in distribution but the limit process $C$ is not of the type of $\mathbb{G}^F$. Denote

$$q_k(t) = I_{\{X_k \le t\}} - k E[I_{\{X_{k+1} \le t\}} | \mathcal{G}_k] + (k-1) E[I_{\{X_k \le t\}} | \mathcal{G}_{k-1}] \qquad \text{for } t \in \mathbb{R}.$$

THEOREM 4.3. *Suppose $(X_n)_n$ is $\mathcal{G}$-c.i.d., $(C_n)_n$ meets condition* (17) *and* $\sup_n E[C_{n,t}^2] < \infty$ *for all $t \in \mathbb{R}$. If*

$$(18) \qquad \frac{1}{n} \sum_{k=1}^{n} k^2 (E[I_{\{X_{k+1} \le t\}} | \mathcal{G}_k] - E[I_{\{X_k \le t\}} | \mathcal{G}_{k-1}])^2 \to 0$$

*in probability $\forall \, t \in \mathbb{R}$,*

*then $C_n \to \mathbb{G}^F$ in distribution and $\mathbb{G}^F$ is tight. Moreover, if*

$$(19) \qquad \frac{1}{n} \sum_{k=1}^{n} q_k(s) q_k(t) \to \sigma(s, t) \qquad \text{a.s. for all } s, t \in \mathbb{R},$$

*then $C_n \to C$ in distribution, where $C$ is a tight process whose distribution $\nu$ is given by*

$$\nu\{x \in \mathcal{X} : (x(t_1), \dots, x(t_r)) \in A\} = \int \mathcal{N}(0, \Sigma(t_1, \dots, t_r))(A) \, dP$$

*for all $t_1, \dots, t_r \in \mathbb{R}$ and $A \in \mathcal{B}(\mathbb{R}^r)$, with $\Sigma(t_1, \dots, t_r) = (\sigma(t_i, t_j))_{1 \le i, j \le r}$.*

PROOF. Arguing as in the proof of Theorem 4.2, it is enough to see that the finite-dimensional distributions of $(C_n)_n$ converge weakly to those of $\mathbb{G}^F$ under (18) and to those of $C$ under (19). Fix $t_1, \dots, t_r, a_1, \dots, a_r \in \mathbb{R}$, define $f = \sum_{i=1}^{r} a_i I_{(-\infty, t_i]}$ and note that $(C_{n,f})_n$ is bounded in $L^2$. If (18) holds, then

$$\frac{1}{n} \sum_{k=1}^{n} k^2 (E[f(X_{k+1})|\mathcal{G}_k] - E[f(X_k)|\mathcal{G}_{k-1}])^2 \to 0 \qquad \text{in probability,}$$

and Theorem 3.3 yields

$$\sum_{i=1}^{r} a_i C_{n,t_i} = C_{n,f} \to \sum_{i=1}^{r} a_i \mathbb{G}_{t_i}^F \qquad \text{in distribution.}$$

Similarly, if (19) holds, then $M_n \to \sum_{i,j} a_i a_j \sigma(t_i, t_j)$ a.s., and Theorem 3.3 implies that $(C_{n,f})_n$ converges in distribution to the probability law $\mu$ on $\mathcal{B}(\mathbb{R})$ given by

$$\mu(A) = \int \mathcal{N}\left(0, \sum_{i,j} a_i a_j \sigma(t_i, t_j)\right)(A)\, dP$$

$$= \nu\left\{x \in \mathcal{X} : \sum_{i=1}^{r} a_i x(t_i) \in A\right\}, \qquad A \in \mathcal{B}(\mathbb{R}).$$

By letting $a_1, \ldots, a_r$ vary, it follows that $(C_{n,t_1}, \ldots, C_{n,t_r}) \to (\mathbb{G}_{t_1}^F, \ldots, \mathbb{G}_{t_r}^F)$ in distribution under (18), and that $(C_{n,t_1}, \ldots, C_{n,t_r}) \to (C_{t_1}, \ldots, C_{t_r})$ in distribution under (19). $\square$

REMARK 4.4. Suppose $(X_n)_n$ is $\mathcal{G}$-c.i.d. and $\mathcal{K} \subset \mathcal{H}$ is a countable class of functions such that $\sup_{f \in \mathcal{K}} |f|$ is in $\mathcal{H}$. Since $C_{n,f} = E[W_{n,f}|\mathcal{G}_n]$ a.s., one obtains

$$E\left[\sup_{f \in \mathcal{K}} |C_{n,f}|\right] \leq E\left[\sup_{f \in \mathcal{K}} |W_{n,f}|\right].$$

Likewise, a direct calculation shows that

$$\limsup_n E\left[\sup_{f \in \mathcal{K}} |B_{n,f}|\right] \leq 5 \limsup_n E\left[\sup_{f \in \mathcal{K}} |W_{n,f}|\right].$$

By these inequalities, $(B_n)_n$ and $(C_n)_n$ can be connected to $(W_n)_n$. In particular, suppose that, for all $\varepsilon > 0$, there is a finite partition $I_1, \ldots, I_m$ of $\mathbb{R}$ by right-open intervals such that

$$(20) \qquad \limsup_n E\left(\max_k \sup_{s,t \in I_k} |W_{n,s} - W_{n,t}|\right) < \varepsilon.$$

Then, (20) still holds with $(B_n)_n$ or $(C_n)_n$ in the place of $(W_n)_n$ and, thus, $(B_n)_n$ and $(C_n)_n$ meet condition (17).

Theorems 4.2 and 4.3 are general results on $\mathcal{G}$-c.i.d. sequences. In the exchangeable case, however, something more can be said. We close the paper by dealing with this case.

THEOREM 4.5.   *If $(X_n)_n$ is exchangeable and $\mathbb{G}^F$ is tight, then $W_n \to \mathbb{G}^F$ in distribution.*

Theorem 4.5 can be proved by a standard application of de Finetti's representation theorem. We refer to Berti, Pratelli and Rigo (2002) for a proof. The assumption that $\mathbb{G}^F$ is tight, instead, needs two remarks. First, it can not be suppressed. Indeed, when $\mathbb{G}^F$ is not tight, $(W_n)_n$ can fail to converge in distribution even if $(X_n)_n$ is exchangeable. An example is in Berti and Rigo (2004). Second, a tight version of $\mathbb{G}^F$ is available if $P(X_1 = X_2) = 0$ or if $P(X_1 \in A) = 1$ for some countable $A \subset \mathbb{R}$. In fact, if $(X_n)_n$ is exchangeable and $P(X_1 = X_2) = 0$, the random distribution function $F$ can be taken to be continuous. Precisely, in some probability space, there are a standard Brownian bridge $\mathbb{G}^0$ and a version of $F$, independent of $\mathbb{G}^0$, whose paths are continuous. Hence, $\mathbb{G}^F = \mathbb{G}^0_F$ is tight. Up to replacing "continuous" with "continuous on $A^c$," the same is true if $P(X_1 \in A) = 1$ for some countable $A$.

Finally, let us turn to $B_n$ and $C_n$. Investigating their asymptotic behavior needs a little more than a straightforward application of de Finetti's theorem.

THEOREM 4.6.   *Suppose $(X_n)_n$ is exchangeable and $\mathbb{G}^F$ is tight. Then, $B_n \to \mathbb{G}^F$ in distribution and $(C_n)_n$ is relatively sequentially compact. Moreover, $C_n \to \mathbb{G}^F$ in distribution under condition (18) and $C_n \to C$ in distribution under condition (19), where $C$ is the tight process described in Theorem 4.3.*

PROOF.   Suppose first that $(W_n)_n$ meets condition (20). Then, by Remark 4.4, $(B_n)_n$ and $(C_n)_n$ satisfy (17). By exchangeability, $(C_{n,t})_n$ is bounded in $L^2$ for all $t$. Thus, Theorems 4.2 and 4.3 yield $B_n \to \mathbb{G}^F$ in distribution, $C_n \to \mathbb{G}^F$ in distribution under (18) and $C_n \to C$ in distribution under (19). Moreover, $(C_n)_n$ is relatively sequentially compact by Lemma 1.5.2 and Theorems 1.3.9 and 1.5.6 of van der Vaart and Wellner (1996). Hence, it is enough to prove (20). If $(X_n)_n$ is i.i.d., then $\sup_n E[\|W_n\|^2] \le c$, where the constant $c$ does not depend on the distribution of $X_1$; see van der Vaart and Wellner (1996), pages 247 and 248. By de Finetti's theorem, $\sup_n E[\|W_n\|^2] \le c$ still holds if $(X_n)_n$ is exchangeable. It follows that $(\|W_n\|)_n$ is uniformly integrable, and since $W_n \to \mathbb{G}^F$ in distribution (by Theorem 4.5), $E[\|\mathbb{G}^F\|] = \lim_n E[\|W_n\|] < \infty$. By Lemma 4.1, given $\varepsilon > 0$, there is a partition $I_1, \dots, I_m$ of $\mathbb{R}$ by right-open intervals such that $E(\max_k \sup_{s,t \in I_k} |\mathbb{G}^F_s - \mathbb{G}^F_t|) < \varepsilon$. Let $h(x) = \max_k \sup_{s,t \in I_k} |x(s) - x(t)|$ for $x \in \mathcal{X}$. Since $h$ is continuous, $h(W_n) \to h(\mathbb{G}^F)$ in distribution. Since $h(W_n) \le 2\|W_n\|$, $(h(W_n))_n$ is uniformly integrable. Thus, $\limsup_n E[h(W_n)] = E[h(\mathbb{G}^F)] < \varepsilon$. $\quad\square$

P. BERTI
DIPARTIMENTO DI MATEMATICA
  PURA ED APPLICATA "G. VITALI"
UNIVERSITA' DI MODENA
  E REGGIO-EMILIA
VIA CAMPI 213/B
41100 MODENA
ITALY
E-MAIL: berti.patrizia@unimore.it

L. PRATELLI
ACCADEMIA NAVALE
VIALE ITALIA 72
57100 LIVORNO
ITALY
E-MAIL: pratel@mail.dm.unipi.it

     P. BERTI, L. PRATELLI AND P. RIGO

P. Rigo
Dipartimento di Economia
    Politica e Metodi Quantitativi
Universita' di Pavia
via S. Felice 5
27100 Pavia
Italy
e-mail: prigo@eco.unipv.it